\documentclass[11pt]{article}
\usepackage{amssymb,amsmath,latexsym, verbatim,enumerate}
\usepackage[hyperfootnotes=false]{hyperref}

\usepackage{tikz}
\usepackage{hyperref}
\allowdisplaybreaks

\begin{document}

\begin{doublespace}

\newtheorem{thm}{Theorem}[section]
\newtheorem{lemma}[thm]{Lemma}
\newtheorem{cond}[thm]{Condition}
\newtheorem{defn}[thm]{Definition}
\newtheorem{prop}[thm]{Proposition}
\newtheorem{proposition}[thm]{Proposition}
\newtheorem{corollary}[thm]{Corollary}
\newtheorem{remark}[thm]{Remark}
\newtheorem{example}[thm]{Example}
\newtheorem{conj}[thm]{Conjecture}
\numberwithin{equation}{section}
\def\ee{\varepsilon}
\def\qed{{\hfill $\Box$ \bigskip}}
\def\NN{{\cal N}}
\def\AA{{\cal A}}
\def\MM{{\cal M}}
\def\BB{{\cal B}}
\def\CC{{\cal C}}
\def\LL{{\cal L}}
\def\DD{{\cal D}}
\def\FF{{\cal F}}
\def\EE{{\cal E}}
\def\QQ{{\cal Q}}
\def\RR{{\mathcal R}}
\def\R{{\mathbb R}}
\def\L{{\mathcal L}}
\def\K{{\bf K}}
\def\S{{\bf S}}
\def\A{{\bf A}}
\def\E{{\mathbb E}}
\def\F{{\bf F}}
\def\P{{\mathbb P}}
\def\N{{\mathbb N}}
\def\eps{\varepsilon}
\def\wh{\widehat}
\def\wt{\widetilde}
\def\pf{\noindent{\bf Proof.} }
\def\beq{\begin{equation}}
\def\eeq{\end{equation}}
\def\lam{\lambda}
\def\H{\mathcal{H}}
\def\nn{\nonumber}
\def\C{\mathbb{C}}
\def\la{\langle}
\def\ra{\rangle}
\def\ud{\mathrm{d}}
\def\hp{\textcolor{blue}}

\newcommand\blfootnote[1]{%
  \begingroup
  \renewcommand\thefootnote{}\footnote{#1}%
  \addtocounter{footnote}{-1}%
  \endgroup
}

\title{\Large \bf  A unified approach to the small-time behavior of the spectral heat content for isotropic L\'evy processes}
\author{Kei Kobayashi and Hyunchul Park}

\date{\today}
\maketitle


\begin{abstract}
This paper establishes the precise small-time asymptotic behavior of the spectral heat content for isotropic L\'evy processes on bounded $C^{1,1}$ open sets of $\R^{d}$ with $d\ge 2$, where the underlying characteristic exponents are regularly varying at infinity with index $\alpha\in (1,2]$, including the case $\alpha=2$. 
Moreover, this asymptotic behavior is shown to be stable under an integrable perturbation of its L\'evy measure. 
These results cover a wide class of isotropic L\'evy processes, including Brownian motions, stable processes, and jump diffusions, and the proofs provide a unified approach to the asymptotic behavior of the spectral heat content for all of these processes. 
\end{abstract}

\section{Introduction}\label{section:introduction}

The spectral heat content (SHC) measures the total heat contained in a domain $D$ at a given time $t>0$ when the initial temperature is set to be one and the temperature outside $D$ is kept at zero.
The SHC contains geometric information of the domain and spectral information of the infinitesimal generator of the underlying process. 
The SHC for Brownian motions was studied intensively more than three decades ago, while the investigation into the SHC for jump processes started within the last decade. In \cite{Val2016, Val2017}, the author established the small-time asymptotic behavior of the SHC for isotropic stable processes on bounded $C^{1,1}$ open sets and provided two-term upper and lower bounds. Later in \cite{PS22}, the precise two-term asymptotic expansion of the SHC for the same processes was established. 
Various questions about the SHC, including those for other jump processes or for higher order terms, have been answered in \cite{GPS19, KP22, KP23, P21, PS19}. 

The investigations to be conducted in this paper are natural continuation of those in \cite{GPS19, PS22}. In \cite[Theorem 4.2]{GPS19}, the authors established the two-term small-time asymptotic expansion of the SHC for L\'evy processes on $\R$ whose characteristic exponents are regularly varying with index $\alpha\in (1,2]$, whereas in \cite[Theorem 1.1]{PS22}, the same problem was answered for isotropic (rotationally invariant) stable processes on $\R^{d}$ with $d\ge 2$ with characteristic exponents $|\xi|^{\alpha}$, $\alpha\in (0,2)$. On the other hand, 
Theorem \ref{thm:main1} of the current paper establishes the two-term small-time asymptotic expansion of the SHC for isotropic L\'evy processes on $\R^{d}$ whose characteristic exponents are regularly varying with index $\alpha\in (1,2]$, thus extending \cite[Theorem 4.2]{GPS19} to higher dimensions while generalizing \cite[Theorem 1.1]{PS22} to cover a large class of L\'evy processes, including non-stable ones. Moreover, Corollary \ref{cor:main} establishes that this asymptotic expansion is stable when the L\'evy measure is perturbed so that the difference of the original L\'evy measure and the new one is integrable. This makes our asymptotic result applicable to the so-called \textit{truncated jump processes}, where all large jumps (say, jump sizes bigger than 1) are removed.

Let us clarify the contributions of this paper to the existing literature. First, the main theorem, Theorem \ref{thm:main1}, is proved under minimal assumptions. We only require that the L\'evy process be isotropic (rotationally invariant) and its characteristic exponent be regularly varying at infinity with index $\alpha\in (1,2]$; no additional assumptions are imposed on the exponent or the associated L\'evy measure. Thus, it covers a large class of L\'evy processes including Brownian motions, stable processes, relativistic stable processes (also called tempered stable processes), mixed stable processes, jump diffusions, and their truncated versions. A list of concrete examples is provided in Section \ref{section:examples}. 
Second, the proof of Theorem \ref{thm:main1} is robust. In fact, when the two-term asymptotic behavior of the SHC for Brownian motions was obtained in \cite[Theorem 6.2]{vdBD89}, the independence of the components played a critical role in their proof, but the latter approach does not work for jump processes since the components are no longer independent. Instead, we will build our discussion upon the methods used to derive \cite[Theorem 1.1]{PS22} for stable processes, while significantly modifying and sophisticating their arguments to overcome the challenges caused by the lack of self-similarity for non-stable L\'evy processes. Consequently, the methods to be employed in this paper can be regarded as an alternative, intricate approach to deriving both \cite[Theorem 6.2]{vdBD89} and \cite[Theorem 1.1]{PS22}. 

Now we briefly illustrate our strategy to prove Theorem \ref{thm:main1}, which indicates that the \textit{heat loss} $|D|-Q_{D}^{X}(t)$ 
(representing the total amount of heat that has been lost from the domain $D$ by time $t>0$ due to the zero exterior condition) decays at the rate of $\psi^{-1}(1/t)^{-1}$, where $\psi$ denotes the regularly varying characteristic exponent of the underlying L\'evy process $X=\{X_t\}_{t\ge 0}$. 
We first verify in Lemma \ref{lemma:inside} that for heat particles starting at points well inside the domain $D$, the heat loss is of order $O(t)$, which is negligible when compared to the above-mentioned rate function $\psi^{-1}(1/t)^{-1}$.  
Next, for heat particles starting at points near the boundary of the domain, we analyze the heat loss by utilizing the \textit{scaled process} $\widetilde{X}^{(t)}=\{\widetilde{X}^{(t)}_u\}_{u\ge 0}$ defined by $\widetilde{X}_{u}^{(t)}=\psi^{-1}(1/t)X_{tu}$ (see Equation \eqref{def:Xtilde}).
For stable processes, this scaled process is equal in law to the original stable process due to self-similarity, so computations are often carried out in a straightforward manner. In our case, however, this is no longer the case, and as a result, various aspects of our discussion not only become complicated but also require new estimates that were not present in \cite{PS22} in the case of stable processes. 
The sharp asymptotic upper bound of the heat loss is established via Lemmas \ref{lemma:half-space}, \ref{lemma:ball}, and \ref{lemma:cancellation}, while the sharp asymptotic lower bound is derived via Lemmas \ref{lemma:half-space}, \ref{lemma:outer ball}, and \ref{lemma:cancellation2}. 
For the upper bound, we use the \textit{interior ball condition} of the $C^{1,1}$ open set $D$; for any $x\in D$ close enough to the boundary, there exists a ball of a fixed radius $R>0$ that contains $x$ and is contained in $D$. 
This property, together with Lemmas \ref{lemma:half-space} and \ref{lemma:ball}, allows us to translate the problem about the heat loss from the domain $D$ into the problem about the heat loss from the upper half-space $H=\{x=(x_1,\ldots, x_d)\in \R^{d}: x_{d}>0\}$.
For the lower bound, in contrast, we use the \textit{exterior ball condition}; for any $x\in D$ close enough to the boundary, there exists a ball of a fixed radius $R>0$ that is contained in $\R^d\setminus\overline{D}$.
Using Lemmas \ref{lemma:outer ball} and \ref{lemma:cancellation2}, we again analyze the heat loss from $D$ by connecting it with the heat loss from the half-space. 

The organization of this paper is as follows. Section \ref{section:preliminaries} sets up notations and recalls necessary facts. The main results, Theorem \ref{thm:main1} and Corollary \ref{cor:main}, are presented and proved in Section \ref{section:main}. Section \ref{section:examples} illustrates concrete examples of L\'evy processes to which our results are applicable.  
The notation $\P_x$ stands for the law of the underlying process started at $x \in \R^{d}$ or $x \in \R$ depending on the context, and $\E_x$ stands for expectation under $\P_x$.
For simplicity, we write $\P = \P_0$ and $\E = \E_0$.

\section{Preliminaries}\label{section:preliminaries}
Let $X=\{X_{t}\}_{t\geq 0}$ be a L\'evy process in $\R^d$ with $d\ge 2$.
By the L\'evy-Khintchine theorem (\cite[Theorem 8.1]{Sato}), 
\[
\E[e^{i\xi \cdot X_{t}}]=e^{-t\psi(\xi)}, \quad \xi\in\R^{d},
\]
with
\[
\psi(\xi)=-i \xi\cdot \gamma + \xi \cdot A\xi -\int_{\R^{d}}\left(e^{i \xi \cdot x}-1-i \xi \cdot x I_{|x|\le 1} \right)\nu(\ud x),
\]
where $\gamma\in\R^{d}$, $A$ is a symmetric, nonnegative-definite matrix, and $\nu$ is a measure satisfying $\nu(\{0\})=0$ and $\int_{\R^d}(1\wedge |x|^2)\nu(\ud x)<\infty$.
The law of $X$ is uniquely characterized by the characteristic exponent $\psi(\xi)$, and the triplet $(A,\gamma,\nu)$ is called the L\'evy triplet of $X$. 
Throughout the paper, $X$ is assumed to be isotropic (i.e.,\ rotationally invariant), which implies that $A=aI$, a constant multiple of the identity matrix, and that $\psi$ is a radial function, so we will write $\psi(b)=\psi(\xi)$ for any $\xi\in \R^{d}$ with $|\xi|=b$.

An isotropic stable L\'evy process with stability index $\alpha$ is a L\'evy process with $A=0$ and $\psi(\xi)=|\xi|^{\alpha}$. 
It is well-known that $\alpha\in (0,2]$, and when $\alpha=2$, it is a Brownian motion that runs twice as fast as the standard Brownian motion. 
In this paper, an isotropic stable process with stability index $\alpha$ in $\R^d$ with $d\ge 2$ is denoted by $X^{(\alpha)}=\{X_{t}^{(\alpha)}\}_{t\geq 0}$, whereas a symmetric stable process in $\R$ is denoted by $Y^{(\alpha)}=\{Y_{t}^{(\alpha)}\}_{t\geq 0}$. 

A real valued function $f$ is said to be regularly varying at infinity with index $\alpha$ if 
\[
\lim_{y\to\infty}\frac{f(xy)}{f(y)}=x^{\alpha} \text{ for all } x>0.
\]
The class of regularly varying functions at infinity with index $\alpha$ is denoted by $\RR_{\alpha}(\infty)$.
For the rest of this paper, we assume that $X$ is an isotropic L\'evy process in $\R^{d}$ with the characteristic exponent $\psi(\xi)\in \RR_{\alpha}(\infty)$ with some $\alpha\in (1,2]$.
Notice that the asymptotic behavior of the spectral heat content when $\alpha\in (0,1)$ was already studied in \cite{GPS19}, and hence, we are left with the case $\alpha\in [1,2]$. In the present paper, however, we restrict our discussion to $\alpha\in (1,2]$ since Theorem \ref{thm:main1} requires the expectation appearing in \eqref{eqn:finite mean} to be finite. On the other hand, the threshold case $\alpha=1$ is worth investigating, and this will be explored in a future project.

Let $D$ be an open set in $\R^{d}$, and let $\tau^X_{D}=\inf\{t>0:X_{t}\notin D\}$ be the first exit time of an isotropic L\'evy process $X$ from $D$. The spectral heat content $Q^X_{D}(t)$ of $X$ at time $t>0$ is defined by 
\[
Q^X_{D}(t)=\int_{D}\P_{x}(\tau_{D}^X>t)\ud x.
\]
Note that $\P_{x}(\tau^X_{D}>t)$ is the unique solution to the following Cauchy problem:
\[
\begin{cases}
\partial_{t}u(t,x)=\mathcal{L}u(t,x),  & (t,x)\in (0,\infty)\times D,\\
\lim_{t\to 0}u(t,x)=1 & x\in D,\\
u(t,y)=0, & t>0,\, y\in D^{c},
\end{cases}
\]
where $\mathcal{L}$ is the infinitesimal generator of the killed process $X^{D}$. Hence, the spectral heat content $Q^X_{D}(t)$ measures the total heat that remains in $D$ at time $t>0$ when the initial temperature of $D$ is set at one and the temperature of $D^{c}$ is kept identically zero.

An open set $D$ in $\R^d$ is said to be $C^{1, 1}$ if 
there exist a localization radius $R>0$ and a constant $\Lambda>0$ that satisfy the following condition (see e.g., \cite{PS22}):
for every $z\in \partial D$, there exist a $C^{1, 1}$ function $\phi=\phi_z: \R^{d-1}\to \R$ satisfying $\phi(0)=0$, $\nabla \phi(0)=(0, \ldots, 0)$,
$\|\nabla\phi\|_\infty\le \Lambda$,
$|\nabla \phi(x_1)-\nabla \phi(x_2)|\le \Lambda |x_1-x_2|$ for $x_1,x_2\in \R^{d-1}$, and an orthonormal coordinate system $CS_z: y=(y_1,\ldots,y_{d-1},y_d)=(\tilde{y}, y_d)$
with origin at $z$ such that
\[
B(z, R)\cap D=B(z, R)\cap \{ y=(\tilde{y}, y_d) \mbox{ in } CS_z: y_d>\phi(\tilde{y})\}.
\]
The pair $(R, \Lambda)$ is called the $C^{1,1}$ characteristics of the $C^{1, 1}$ open set $D$. 
It is well-known that any $C^{1, 1}$ open set $D$ with $C^{1,1}$ characteristics $(R,\Lambda)$ satisfies the uniform interior and exterior $R$-ball condition 
(see \cite[Lemma 2.2]{AKSZ2007}): 
for any $z\in\partial D$,  there exist open balls $B_1$ and $B_2$ of the same radius $R$ such that 
\begin{align}\label{def:interior-exterior}
	B_1\subset D, \ B_2\subset\R^d\setminus\overline{D}, \ \textrm{and} \ \partial B_1\cap\partial B_2=\{z\}.
\end{align}

For $x\in D$, let $\delta_{D}(x)$ denote the Euclidean distance from $x$ to the boundary $\partial D$; i.e.,
\begin{align}\label{def:deltaD}
	\delta_{D}(x)=\inf\{|x-y| : y \in \partial D\}. 
\end{align}
For each $a\in(0,R/2]$, let
\begin{align}\label{def:Da}
	D_{a}=\{x\in D : \delta_{D}(x)>a\},
\end{align}
which is a region obtained by removing points near the boundary $\partial D$ from $D$. 
Note that the $C^{1,1}$ condition with $C^{1,1}$ characteristics $(R,\Lambda)$ implies that the set $D_a$ is non-empty. 
It also follows from \cite[Lemma 6.7]{vdBD89} that if the $C^{1,1}$ open set $D$ is bounded, then
\beq\label{eqn:vdBD89}
|\partial D|\left(\frac{R-a}{R}\right)^{d-1}\leq |\partial D_{a}|\leq |\partial D|\left(\frac{R}{R-a}\right)^{d-1}, 
\eeq
where $|\partial D|$ and $|\partial D_a|$ denote the perimeters (or the Lebesgue measures) of the sets $\partial D$ and $\partial D_a$, respectively. 
In particular, $|\partial D_{a}|\leq 2^{d-1}|\partial D|$ for any $a\in(0,R/2]$. 
The reader is advised to note that the same notation $|\cdot|$ is used throughout the paper to denote the Euclidean norm.

For $x\in D\setminus D_{R/2}$ (so that $\delta_D(x)\le R/2$), let $z_x$ be the unique point in $\partial D$ such that $|x-z_x|=\delta_{D}(x)$.
Let $\mathbf{n}_{z_{x}}=\frac{z_{x}-x}{|z_x-x|}$ be the outward unit normal vector to $\partial D$ at the point $z_x$. For an isotropic $\alpha$-stable L\'evy process $X^{(\alpha)}$ starting at $x$, let 
\[
	Y^{(\alpha)}_{t}= (X^{(\alpha)}_{t}-x)\cdot \mathbf{n}_{z_{x}}
\]
be the one-dimensional projection of $X^{(\alpha)}_{t}-x$ onto the direction of $\mathbf{n}_{z_{x}}$. 
Then the process $Y^{(\alpha)}=\{Y^{(\alpha)}_t\}_{t\ge 0}$ is a symmetric $\alpha$-stable process in $\R$. Moreover, in the particular case when $X^{(\alpha)}$ starts at a point $(\tilde{0}, v)=(0,\ldots,0,v)$ in the upper half-space $H=\{x=(\tilde{x}, x_{d})=(x_1,\ldots,x_{d-1},x_d)\in \R^d : x_{d}>0\}$, the relationship 
\begin{align}\label{Xalpha-Yalpha}
	\P_{(\tilde{0},v)}(\tau_{H}^{X^{(\alpha)}}\le 1)=\P(\overline{Y}_{1}^{(\alpha)}\ge v) 
\end{align}
holds, where $\P$ without the subscript means that the process starts from the origin, and
\[
\overline{Y}^{(\alpha)}_{t}:=\sup_{0\leq s\leq t}Y^{(\alpha)}_{s} 
\]
 represents the supremum process of $Y^{(\alpha)}$; see the proofs of \cite[Lemma 3.1, Proposition 3.3]{PS22} for details. By \cite[Proposition 4 in Section VIII.1]{Ber}, there exists a constant $C_\alpha>0$ such that
\[
\P(\overline{Y}^{(\alpha)}_1>u)\sim\P(Y^{(\alpha)}_1>u)\sim C_\alpha u^{-\alpha} 
\text{ as } u\to \infty.
\]
The identity $\E[\overline{Y}_{1}^{(\alpha)}]=\int_{0}^{\infty}\P(\overline{Y}_{1}^{(\alpha)}>u)\ud u$ implies that
$\E[\overline{Y}_{1}^{(\alpha)}] <\infty$ if $\alpha\in (1,2),$ and $\E[\overline{Y}_{1}^{(\alpha)}] =\infty$ if $\alpha\in (0,1].$
The case $\alpha=2$ corresponds to Brownian motion, and it is well-known that $\E[\overline{Y}_{1}^{(2)}]=\E[|Y_1^{(2)}|]=2\int_{0}^{\infty}\frac{x}{\sqrt{4\pi}}e^{-\frac{x^2}{4}}\ud x=\frac{2}{\sqrt{\pi}}<\infty$.
Hence, 
\beq\label{eqn:finite mean}
\E[\overline{Y}_{1}^{(\alpha)}]<\infty \text{ if and only if } \alpha\in (1,2].
\eeq

\section{Main Results}\label{section:main}

The main result of this paper is the following:

\begin{thm}\label{thm:main1}
Let $D$ be a bounded $C^{1,1}$ open set in $\R^{d}$ with $d\ge 2$.
Let $X$ be an isotropic L\'evy process in $\R^d$ with the characteristic exponent $\psi\in \RR_{\alpha}(\infty)$, $\alpha\in (1,2]$. 
Then 
\[
\lim_{t\to0}\psi^{-1}(1/t)(|D|-Q_{D}^{X}(t))=|\partial D|\E[\overline{Y}^{(\alpha)}_{1}],
\]
where $\overline{Y}^{(\alpha)}$ is the supremum process of a symmetric $\alpha$-stable process $Y^{(\alpha)}$ in $\R$.
\end{thm}

Theorem \ref{thm:main1} continues to hold for a L\'evy process $U$ whose L\'evy measure is close to that of $X$ in the sense that the difference of the two is a finite signed measure.
In this case, one can use \cite[Theorem 5.1]{GPS19} whose main idea is to decompose the process $U$ as an independent sum of the process $X$ and a compound Poisson process. Here is a precise statement:
\begin{corollary}\label{cor:main}
Let $D$ be a bounded $C^{1,1}$ open set in $\R^{d}$ with $d\ge 2$.
Let $U$ be an isotropic L\'evy process in $\R^d$ with L\'evy triplet $(aI,0,\nu^{U})$. If there exists another isotropic L\'evy process with L\'evy triplet $(aI,0,\nu^X)$ and the characteristic exponent $\psi\in\RR_{\alpha}(\infty)$, $\alpha\in (1,2]$, such that $\sigma(\ud x)=\nu^U(\ud x)-\nu^X(\ud x)$ is a finite signed measure, then 
\[
|D|-Q_{D}^{U}(t)=
|\partial D|\E[\overline{Y}^{(\alpha)}_{1}]\psi^{-1}(1/t)^{-1} + o(\psi^{-1}(1/t)^{-1}).
\]
\end{corollary}
\pf
The condition $\psi\in \mathcal{R}_{\alpha}(\infty)$  implies $\psi^{-1}\in \mathcal{R}_{1/\alpha}(\infty)$ due to \cite[Theorems 1.5.3 and 1.5.12]{BGT}, and in particular,  
$t=o(\psi^{-1}(1/t)^{-1})$ as $t\to 0$.
Hence, the desired conclusion follows from \cite[Theorem 5.1]{GPS19} and Theorem \ref{thm:main1}.
\qed

The remainder of the section is devoted to the proof of Theorem \ref{thm:main1}. 
Our discussion builds upon the approach taken in \cite{PS22} for stable processes, but we will modify that approach substantially to cover more general isotropic L\'evy processes.

The first lemma below shows that the contribution of points inside the domain $D$ to the amount of heat loss $|D|-Q_{D}^{X}(t)=\int_D \P_x(\tau_{D}^{X}\leq t)\ud x$ is at most of order $t$. 
On the other hand, since $\psi^{-1}\in \mathcal{R}_{1/\alpha}(\infty)$ when $\psi\in \mathcal{R}_{\alpha}(\infty)$ due to \cite[Theorems 1.5.3 and 1.5.12]{BGT}, it follows that 
$t=o(\psi^{-1}(1/t)^{-1})$ 
as $t\to 0$. 
Hence, it is concluded that the contribution to the heat loss from inside the domain $D$ is negligible when it is compared to $\psi^{-1}(1/t)^{-1}$. The lemma will be used later in the derivation of the upper bound for Theorem \ref{thm:main1}. 
Recall the notation $D_a$ introduced in \eqref{def:Da}.

\begin{lemma}\label{lemma:inside}
Let $D$ and $X$ be as in Theorem \ref{thm:main1}, where the $C^{1,1}$ characteristics of $D$ are given by $(R,\Lambda)$. Then for any $a\in(0,R/2]$, there exists a constant $c>0$ such that  
\[
\int_{D_a}\P_{x}(\tau_{D}^{X}\leq t)\ud x \leq c|D|t
\]
for all $t>0$.
\end{lemma}
\pf 
For any $x\in D_a$, by \cite[Equation (3.2)]{Pruitt81}, 
\begin{align*}
\P_{x}( \tau_{D}^X \leq t) \leq \P_x\left(\sup _{s\leq t}|X_{s}-x| \geq a\right)
=\P\left(\sup _{s\leq t}|X_{s}| \geq a\right)
\le c t,
\end{align*}
where $c>0$ is a constant independent of $x$. 
In fact, \cite[Equation (3.2)]{Pruitt81} assumes that the Gaussian part of $X$ vanishes, but one can include a Gaussian part given by a constant multiple of a Brownian motion (recall that $X$ is assumed to be isotropic) since the tail of its supremum 
process has an upper bound that decays exponentially as $t$ decreases; see \cite[Chapter 2, Equation (8.3)']{KaratzasShreve}. 
Therefore,
\[
\int_{D_{a}}\P_{x}(\tau_{D}^{X}\leq t)\ud x\leq c|D_{a}|t \leq c|D|t,
\]
which completes the proof.
\qed

The next lemma is similar to \cite[Proposition 3.3]{PS22}.
Let 
\[
	H=\{x=(x_{1},\ldots, x_{d})\in \R^d : x_{d}>0\}
\]
 be the upper half-space.
 In the remainder of the paper, we use the notation $\tilde{0}$ to represent the zero vector in $\R^{d-1}$.

\begin{lemma}\label{lemma:half-space}
Let $D$ and $X$ be as in Theorem \ref{thm:main1}, where the $C^{1,1}$ characteristics of $D$ are given by $(R,\Lambda)$.  
Then there exists $A=A(\psi)\in(0,R/2]$ such that for any $a\in(0, A ]$,
\[
\lim_{t\to 0}\psi^{-1}(1/t)\int_{0}^{a}\P_{(\tilde{0},u)}(\tau^{X}_{H} \leq t)\ud u=\E[\overline{Y}^{(\alpha)}_{1}].
\]
\end{lemma}
\pf
For each $t>0$, define the scaled process $\widetilde{X}^{(t)}=\{\widetilde{X}^{(t)}_{u}\}_{u\geq 0}$ by
\begin{align}\label{def:Xtilde}
\widetilde{X}_{u}^{(t)}=\psi^{-1}(1/t)X_{tu}.
\end{align}
Referring to the proof of \cite[Lemma 4.5]{GPS19}, one can verify that the characteristic exponent $t\psi(\psi^{-1}(1/t)\xi)$ of $\widetilde{X}^{(t)}$ converges to the characteristic exponent $|\xi|^\alpha$ of an isotropic $\alpha$-stable process $X^{(\alpha)}$ as $t\to 0$. By \cite[Exercise 16.10]{Kallenberg2002}, the latter implies that $\widetilde{X}^{(t)}$ converges weakly to $X^{(\alpha)}$ in the Skorokhod space equipped with the $J_1$-topology.
Note that 
\begin{align*}
\tau_{H}^{X}
&=\inf\{u : X_{u}\notin H\}
=\inf\{tu: \psi^{-1}(1/t) X_{tu}\notin \psi^{-1}(1/t) H\}\\
&=\inf\{tu: \psi^{-1}(1/t)X_{tu}\notin H\}
=t\inf\{u: \widetilde{X}_{u}^{(t)} \notin H\}=t\tau_{H}^{\widetilde{X}^{(t)}},
\end{align*}
where we used the fact that the upper half-space $H$ is invariant under the multiplication by any positive constant. The latter
 shows that the law of $\tau_{H}^{X}$ under $\P_{x}$ is equal to the law of $t\tau_{H}^{\widetilde{X}^{(t)}}$ under $\P_{\psi^{-1}(1/t)x}$.
Hence, for a fixed $a\in(0,R/2]$, by the change of variables $v=\psi^{-1}(1/t)u$,
\begin{align}\label{eqn:half-space}
\int_{0}^{a}\P_{(\tilde{0},u)}(\tau_{H}^X \leq t)\ud u
&=\int_{0}^{a}\P_{(\tilde{0},\psi^{-1}(1/t)u)}(\tau_{H}^{\widetilde{X}^{(t)}} \leq 1)\ud u\nn\\
&=\psi^{-1}(1/t)^{-1}\int_{0}^{a\psi^{-1}(1/t)}\P_{(\tilde{0},v)}(\tau_{H}^{\widetilde{X}^{(t)}}\leq  1)\ud v.
\end{align}
By \eqref{eqn:finite mean}, $\E[\overline{Y}_{1}^{(\alpha)}]=\int_{0}^{\infty}\P(\overline{Y}_{1}^{(\alpha)}\geq v)\ud v<\infty$ for $\alpha\in (1,2]$. Hence, for any $\eps>0$, there exists $N=N(\eps)>0$ such that 
\[
\int_{0}^{N}\P(\overline{Y}_{1}^{(\alpha)}\geq v)\ud v>\E[\overline{Y}_{1}^{(\alpha)}]-\eps.
\]
Let 
\[
\widetilde{Z}_{u}^{(t)}=\widetilde{X}_{u}^{(t)}\cdot e_{d}=\psi^{-1}(1/t)X_{tu}\cdot e_{d}
\]
be the projection of $\widetilde{X}_{u}^{(t)}$ onto the last coordinate axis. The characteristic function of $\widetilde{Z}_{1}^{(t)}$ is
\[
	\E[e^{i\xi\widetilde{Z}_{1}^{(t)}}]=\E[e^{i\xi \psi^{-1}(1/t)X_{t}\cdot e_{d}}]=e^{-t\psi(\psi^{-1}(1/t)\xi\cdot e_{d})}
	=e^{-t\psi(\psi^{-1}(1/t)\xi)},
\]
where the reader is reminded that $\psi$ is a radial function and that we write $\psi(b)=\psi(\xi)$ for any $\xi\in \R^{d}$ with $|\xi|=b$; the latter implies weak convergence of $\widetilde{Z}^{(t)}$ to $Y^{(\alpha)}$ as $t\to 0$ with respect to the $J_1$-topology. 
Note also that since $\psi^{-1}(1/t)\to \infty$ as $t\to 0$, there exists $t_{0}>0$ such that $a\psi^{-1}(1/t) \geq N$ for all $0<t\leq t_{0}$.
It now follows from \eqref{eqn:half-space}, Fatou's lemma, and the weak convergence of $\widetilde{Z}^{(t)}$ to $Y^{(\alpha)}$ that
\begin{align*}
\liminf_{t\to 0}\psi^{-1}(1/t)\int_{0}^{a}\P_{(\tilde{0},u)}(\tau_{H}^X\leq t)\ud u
&\geq \liminf_{t\to 0}\int_{0}^{N}\P_{(\tilde{0},v)}(\tau_{H}^{\widetilde{X}^{(t)}}\leq 1)\ud v\\
&\ge \liminf_{t\to 0}\int_{0}^{N}\P(\overline{\widetilde{Z}^{(t)}}_1 > v)\ud v\\
&\ge \int_{0}^{N}\P(\overline{Y}_{1}^{(\alpha)}> v)\ud v 
>\E[\overline{Y}_{1}^{(\alpha)}]-\eps,
\end{align*}
where $\overline{\widetilde{Z}^{(t)}}_1=\sup_{0\le u\le 1}\widetilde{Z}^{(t)}_u$.
Since $\eps>0$ is arbitrary, we obtain the lower bound  
\begin{align*}
\liminf_{t\to 0}\psi^{-1}(1/t)\int_{0}^{a}\P_{(\tilde{0},u)}(\tau_{H}^X\leq t)\ud u\geq \E[\overline{Y}_{1}^{(\alpha)}].
\end{align*}

Derivation of the upper bound requires a delicate discussion. In fact, a simple modification of the above argument  would be to notice by \eqref{eqn:half-space} that 
\[
\limsup_{t\to 0}\psi^{-1}(1/t)\int_{0}^{a}\P_{(\tilde{0},u)}(\tau_{H}^X \leq t)\ud u
= \limsup_{t\to 0}\int_{0}^{\infty}\P_{(\tilde{0},v)}(\tau_{H}^{\widetilde{X}^{(t)}}\le 1)\mathbf{1}_{(0,a\psi^{-1}(1/t)]}(v)\ud v
\]
and use the reverse version of Fatou's lemma for the latter limit. However, that would require the integrand $\P_{(\tilde{0},v)}(\tau_{H}^{\widetilde{X}^{(t)}}\le 1)\mathbf{1}_{(0,a\psi^{-1}(1/t)]}(v)$ to be bounded above by a function of $v$ that is both independent of $t$ and integrable on the unbounded interval $(0,\infty)$, and finding such an upper bound is a non-trivial task. 
To overcome this hurdle, we adopt the approach employed in the proof of \cite[Lemma 4.5]{GPS19} and modify it to suit our context. First, fix $\delta\in(0,\alpha-1)$. Since $\psi\in \RR_{\alpha}(\infty)$, it follows from Potter's theorem (\cite[Theorem 1.5.6]{BGT}) that there exists $x_{0}=x_{0}(\psi)>0$ such that
\[
\frac{\psi(\psi^{-1}(1/t)1/u)}{\psi(\psi^{-1}(1/t))}\leq \frac{2}{u^{\alpha-\delta}}
\]
for all small $t$ for which $\psi^{-1}(1/t)\geq x_0$ and $1\leq u \leq\frac{\psi^{-1}(1/t)}{x_0}$. (Recall that $\psi^{-1}(1/t)\to \infty$ as $t\to 0$ since $\psi^{-1}\in\mathcal{R}_{1/\alpha}(\infty)$.)
Let 
\[
	A:=\frac{1}{x_0}\wedge \frac{R}{2},
\] 
and consider the above-mentioned one-dimensional projection $\widetilde{Z}^{(t)}$ of $\widetilde{X}^{(t)}$.
Taking $M\ge 1$ and applying the discussion given in the proof of \cite[Lemma 4.5]{GPS19}, one can find a constant $c=c(x_0)>0$ such that for any $a\in (0,A]$ and any $v$ and $t$ satisfying $M\leq v\leq a\psi^{-1}(1/t)$, the inequality 
\beq\label{eqn:ub}
\P(\overline{\widetilde{Z}^{(t)}}_1\geq v) \leq \frac{c}{v^{\alpha-\delta}}
\eeq
holds, where $\overline{\widetilde{Z}^{(t)}}_1=\sup_{0\le u\le 1}\widetilde{Z}^{(t)}_u$.
By \eqref{eqn:ub}, for a fixed $a\in (0,A]$ and any small $t$ for which $M<a\psi^{-1}(1/t)$,
\begin{align}\label{eqn:ub2}
\int_{M}^{a\psi^{-1}(1/t)}\P_{(\tilde{0},v)}(\tau_{H}^{\widetilde{X}^{(t)}}\leq  1)\ud v
&=\int_{M}^{a\psi^{-1}(1/t)}\P(\overline{\widetilde{Z}^{(t)}}_1\geq v) \ud v\notag 
\le \int_{M}^{a\psi^{-1}(1/t)}\frac{c}{v^{\alpha-\delta}}\ud v\\
&\leq \int_{M}^{\infty}\frac{c}{v^{\alpha-\delta}}\ud v
=\frac{c}{\alpha-\delta-1}M^{-\alpha+\delta+1}.
\end{align}

Now, using \eqref{eqn:half-space} and \eqref{eqn:ub2}, applying the reverse version of Fatou's lemma with the trivial upper bound $\P(\overline{\widetilde{Z}^{(t)}}_1\geq v)\leq 1$ which is integrable on the bounded interval $(0,M]$, and using the weak convergence of $\widetilde{Z}^{(t)}$ to $Y^{(\alpha)}$, we obtain
\begin{align*}
\limsup_{t\to 0}\psi^{-1}(1/t)\int_{0}^{a}\P_{(\tilde{0},u)}(\tau_{H}^X \leq t)\ud u
&\le \limsup_{t\to 0}\int_0^{M}\P_{(\tilde{0},v)}(\tau_{H}^{\widetilde{X}^{(t)}}\leq  1)\ud v\\ 
&\quad +\limsup_{t\to 0}\int_{M}^{a\psi^{-1}(1/t)}\P_{(\tilde{0},v)}(\tau_{H}^{\widetilde{X}^{(t)}}\leq  1)\ud v\\
&\le\limsup_{t\to 0} \int_{0}^{M}\P(\overline{\widetilde{Z}^{(t)}}_1\geq v)\ud v +\frac{c}{\alpha-\delta-1}M^{-\alpha+\delta+1}\\
&\le\int_{0}^{M}\P(\overline{Y}_{1}^{(\alpha)}\geq v)\ud v +\frac{c}{\alpha-\delta-1}M^{-\alpha+\delta+1}.
\end{align*}
Letting $M\to \infty$ and using the relationship \eqref{Xalpha-Yalpha} yields the upper bound
\begin{align*}
\limsup_{t\to 0}\psi^{-1}(1/t)\int_{0}^{a}\P_{(\tilde{0},u)}(\tau_{H}^X\leq t)\ud u\leq \E[\overline{Y}_{1}^{(\alpha)}],
\end{align*}
as desired.
\qed

The next lemma shows that one can replace the upper half-space $H$ in Lemma \ref{lemma:half-space} with the ball $B((\tilde{0},R),R)$ of radius $R$ centered at $(\tilde{0},R)\in H$. 
\begin{lemma}\label{lemma:ball}
Let $D$ and $X$ be as in Theorem \ref{thm:main1}, where the $C^{1,1}$ characteristics of $D$ are given by $(R,\Lambda)$. 
Let $A\in(0,R/2]$ be the constant appearing in Lemma \ref{lemma:half-space}.
Then for any $a\in(0, A ]$,
\[
\lim_{t\to 0}\psi^{-1}(1/t)\int_{0}^{a}\P_{(\tilde{0},u)}(\tau^{X}_{B((\tilde{0},R),R)}\leq t)\ud u=\E[\overline{Y}^{(\alpha)}_{1}].
\]
\end{lemma}
\pf
Since $B((\tilde{0},R),R)\subset H$, it follows that $\P_{(\tilde{0},u)}(\tau_{H}^X\leq t) \leq \P_{(\tilde{0},u)}(\tau^{X}_{B((\tilde{0},R),R)}\leq t)$. This, together with Lemma \ref{lemma:half-space}, yields the lower bound
\[
\liminf_{t\to 0}\psi^{-1}(1/t)\int_{0}^{a}\P_{(\tilde{0},u)}(\tau^{X}_{B((\tilde{0},R),R)} \le t)\ud u\geq \E[\overline{Y}^{(\alpha)}_{1}].
\]

To derive the upper bound, recall the definition of $\widetilde{X}^{(t)}$ in \eqref{def:Xtilde} and observe that
\begin{align*}
\tau_{B((\tilde{0},R),R)}^{X}
&=\inf\{u: X_{u}\notin B((\tilde{0},R),R)\}=\inf\{tu: X_{tu}\notin B((\tilde{0},R),R)\}\\
&=\inf\{tu: \psi^{-1}(1/t)X_{tu} \notin \psi^{-1}(1/t)B((\tilde{0},R),R) \}\\
&=t\inf\{u: \widetilde{X}_{u}^{(t)}\notin B((\tilde{0},\psi^{-1}(1/t)R),\psi^{-1}(1/t)R) \}\\
&=t\tau_{B((\tilde{0},\psi^{-1}(1/t)R),\psi^{-1}(1/t)R)}^{\widetilde{X}^{(t)}}.
\end{align*}
Hence, by the change of variables $v=\psi^{-1}(1/t)u$, 
\begin{align}\label{eqn:ball lb1}
\int_{0}^{a}\P_{(\tilde{0},u)}(\tau^{X}_{B((\tilde{0},R),R)}\leq t)\ud u\nn
&=\int_{0}^{a}\P_{(\tilde{0},\psi^{-1}(1/t)u)}(\tau_{B((\tilde{0},\psi^{-1}(1/t)R),\psi^{-1}(1/t)R)}^{\widetilde{X}^{(t)}} \leq 1)\ud u\nn\\
&=\psi^{-1}(1/t)^{-1}\int_{0}^{a\psi^{-1}(1/t)}\P_{(\tilde{0},v)}(\tau_{B((\tilde{0},\psi^{-1}(1/t)R),\psi^{-1}(1/t)R)}^{\widetilde{X}^{(t)}} \leq 1)\ud v.
\end{align}

Note that if $v\in (0, a\psi^{-1}(1/t)]$, then $B((\tilde{0}, v),v)\subset B((\tilde{0},\psi^{-1}(1/t)R),\psi^{-1}(1/t)R)$ since we take $a\le A\le R/2$. 
Moreover, $\max_{1\leq k\leq d}|x_k|\leq |x|=\sqrt{\sum_{k=1}^{d}x_k^2}\leq \sqrt{d}\max_{1\leq k\leq d}|x_k|$ for any $x=(x_1,\ldots, x_d)\in \R^{d}$.
Therefore, for any $v\in (0, a\psi^{-1}(1/t)]$, 
\begin{align}\label{eqn:ball lb4}
\P_{(\tilde{0},v)}(\tau_{B((\tilde{0},\psi^{-1}(1/t)R),\psi^{-1}(1/t)R)}^{\widetilde{X}^{(t)}} \leq 1)
&\leq \P_{(\tilde{0},v)}(\tau_{B((\tilde{0},v),v)}^{\widetilde{X}^{(t)}} \leq 1)\notag 
=\P\left(\sup_{u\leq 1}|\widetilde{X}_{u}^{(t)}|\geq v\right)\\
&\leq\sum_{k=1}^{d}\P\left(\sup_{u\leq 1}|\widetilde{Z}_{u}^{(t),k}|\geq \frac{v}{\sqrt{d}}\right), 
\end{align}
where 
\[
\widetilde{Z}_{u}^{(t),k}=\widetilde{X}_{u}^{(t)}\cdot e_{k}=\psi^{-1}(1/t)X_{tu}\cdot e_{k}
\]
denotes the projection of $\widetilde{X}_{u}^{(t)}$ onto the $k$th coordinate axis.
As in the proof of  Lemma \ref{lemma:half-space}, one can verify that the characteristic function of $\widetilde{Z}_{1}^{(t),k}$ is $\E[e^{i\xi\widetilde{Z}_{1}^{(t),k}}]=e^{-t\psi(\psi^{-1}(1/t)\xi)}$ for all $k$, so the processes $\widetilde{Z}^{(t),k}$, $1\le k\le d$, have the same distribution. 
Hence, the argument given in the proof of Lemma \ref{lemma:half-space} shows that if we fix $\delta\in(0,\alpha-1)$ and $M\ge 1$, then for any $a\in (0,A]$ and any $v$ and $t$ satisfying $M\leq v\leq a\psi^{-1}(1/t)$, the inequality \eqref{eqn:ub} holds with $\widetilde{Z}^{(t),k}$.
This, together with the symmetry of $\widetilde{Z}^{(t),k}$, yields  
\begin{align}\label{eqn:ball lb5}
\sum_{k=1}^{d}\P\left(\sup_{u\leq 1}|\widetilde{Z}_{u}^{(t),k}|\geq \frac{v}{\sqrt{d}}\right)
&\leq d\cdot \P\left(\sup_{u\leq 1}\widetilde{Z}_{u}^{(t),1}\geq \frac{v}{\sqrt{d}} \text{ or } \inf_{u\leq 1}\widetilde{Z}_{u}^{(t),1}\leq -\frac{v}{\sqrt{d}}\right) \notag\\
&=2d\cdot \P\left(\sup_{u\leq 1}\widetilde{Z}_{u}^{(t),1}\geq \frac{v}{\sqrt{d}}\right)
\leq \frac{c}{v^{\alpha-\delta}},
\end{align}
for any $v$ and $t$ satisfying $M \leq \frac{v}{\sqrt{d}}\leq a\psi^{-1}(1/t)$, where $c>0$ is some constant. 
In particular, by \eqref{eqn:ball lb4} and \eqref{eqn:ball lb5}, we deduce that 
\begin{align}\label{eqn:ball lb6}
	\limsup_{t\to0}\int_{M\sqrt{d}}^{a\psi^{-1}(1/t)}\P_{(\tilde{0},v)}(\tau_{B((\tilde{0},\psi^{-1}(1/t)R),\psi^{-1}(1/t)R)}^{\widetilde{X}^{(t)}} \leq 1)\ud v
	&\le \limsup_{t\to0}\int_{M\sqrt{d}}^{a\psi^{-1}(1/t)}\frac{c}{v^{\alpha-\delta}}\ud v \notag\\
	&=\frac{c\bigl(M\sqrt{d}\bigr)^{-\alpha+\delta+1}}{\alpha-\delta-1}.
\end{align}

Now, to deal with the integral $\int_{0}^{M\sqrt{d}}\P_{(\tilde{0},v)}(\tau_{B((\tilde{0},\psi^{-1}(1/t)R),\psi^{-1}(1/t)R)}^{\widetilde{X}^{(t)}} \leq 1)\ud v$, we construct an increasing sequence of domains $D(n)$ by
\[
D(n)=\left\{x=(\tilde{x}, x_{d})\in \R^{d-1}\times \R: |\tilde{x}|<n,\, 0<x_{d}<n,\, \cos^{-1}(\vec{\bf{n}}\cdot \frac{x}{|x|} )<\frac{\pi}{2}-\frac{1}{n} \right\}, \quad \vec{\bf{n}}=(\tilde{0},1).
\]
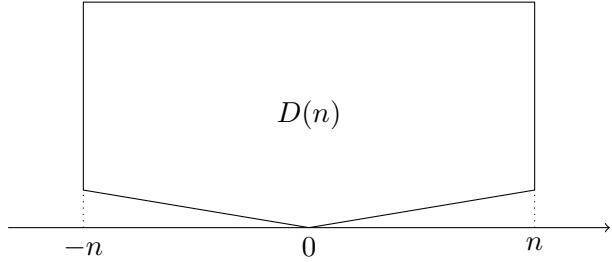
\begin{center}
\begin{tikzpicture}
\draw [->] (-4,0)--(4,0);
\draw (0,0)--(-3,0.5)--(-3,3)--(3,3)--(3,0.5)--(0,0);
\draw [dotted] (3,0)--(3,0.5);
\draw [dotted] (-3,0)--(-3,0.5);
\node [below] at (3,0) {$n$};
\node [below] at (-3,0) {$-n$};
\node [below] at (0,0) {0};
\node at (0,1.5) {$D(n)$};
\end{tikzpicture}
\end{center}
Since $D(n)$ increases to $H$ as $n\to\infty$, it follows from the bounded convergence theorem that
\[
\lim_{n\to\infty}\int_{0}^{M\sqrt{d}}\P_{(\tilde{0},v)}(\tau_{D(n)}^{X^{(\alpha)}}\leq 1)\ud v=\int_{0}^{M\sqrt{d}}\P_{(\tilde{0},v)}(\tau_{H}^{X^{(\alpha)}}\leq 1)\ud v,
\]
where the sequence converges decreasingly.
Therefore, for a given $\eps>0$, we may take an integer $N=N(\eps)$ such that 
\beq\label{eqn:ball lb2}
\int_{0}^{M\sqrt{d}}\P_{(\tilde{0},v)}(\tau_{D(N)}^{X^{(\alpha)}}\leq 1)\ud v<\int_{0}^{M\sqrt{d}}\P_{(\tilde{0},v)}(\tau_{H}^{X^{(\alpha)}}\leq 1)\ud v + \eps.
\eeq
Since $B((\tilde{0},\psi^{-1}(1/t)R),\psi^{-1}(1/t)R)$ increases to $H$ as $t\to 0$, we can take $t_{0}=t_{0}(N)>0$ such that for any $0<t\leq t_{0}$,
\[
D(N)  \subset B((\tilde{0},\psi^{-1}(1/t)R),\psi^{-1}(1/t)R), 
\]
and hence,
\beq\label{eqn:ball lb3}
\P_{(\tilde{0},v)}(\tau_{B((\tilde{0},\psi^{-1}(1/t)R),\psi^{-1}(1/t)R)}^{\widetilde{X}^{(t)}} \leq 1)
\leq \P_{(\tilde{0},v)}(\tau_{D(N)}^{\widetilde{X}^{(t)}}\leq 1).
\eeq
Combining \eqref{eqn:ball lb2} and \eqref{eqn:ball lb3}, applying the reverse version of Fatou's lemma with the trivial upper bound $\P_{(\tilde{0},v)}(\tau_{D(N)}^{\widetilde{X}^{(t)}} \leq 1)\leq 1$ on the bounded interval $(0,M\sqrt{d}]$, and using the weak convergence of $\widetilde{X}^{(t)}$ to $X^{(\alpha)}$, we obtain
\begin{align}\label{eqn:ball lb7}
\limsup_{t\to0}\int_{0}^{M\sqrt{d}}\P_{(\tilde{0},v)}(\tau_{B((\tilde{0},\psi^{-1}(1/t)R),\psi^{-1}(1/t)R)}^{\widetilde{X}^{(t)}} \leq 1)\ud v
&\le\limsup_{t\to0}\int_{0}^{M\sqrt{d}}\P_{(\tilde{0},v)}(\tau_{D(N)}^{\widetilde{X}^{(t)}} \leq 1)\ud v \notag\\ 
&\le\int_{0}^{M\sqrt{d}}\P_{(\tilde{0},v)}(\tau_{D(N)}^{X^{(\alpha)}} \leq 1)\ud v\notag \\
&<\int_{0}^{M\sqrt{d}}\P_{(\tilde{0},v)}(\tau_{H}^{X^{(\alpha)}}\leq 1)\ud v+\eps. 
\end{align}

Finally, putting together \eqref{eqn:ball lb1}, \eqref{eqn:ball lb6} and \eqref{eqn:ball lb7} yields
\begin{align*}
\limsup_{t\to0}\psi^{-1}(1/t)\int_{0}^{a}\P_{(\tilde{0},u)}(\tau^{X}_{B((\tilde{0},R),R)}\leq t)\ud u
<\int_{0}^{M\sqrt{d}}\P_{(\tilde{0},v)}(\tau_{H}^{X^{(\alpha)}}\leq 1)\ud v+\eps+\frac{c\bigl(M\sqrt{d}\bigr)^{-\alpha+\delta+1}}{\alpha-\delta-1}.
\end{align*}
Recalling the relationship \eqref{Xalpha-Yalpha} and letting $\eps\to 0$ and $M\to \infty$ yields
\[
\limsup_{t\to 0}\psi^{-1}(1/t)\int_{0}^{a}\P_{(\tilde{0},u)}(\tau^{X}_{B((\tilde{0},R),R)}\le t)\ud u\leq \E[\overline{Y}^{(\alpha)}_{1}],
\]
which is the desired upper bound.
\qed

Recall the notations $\delta_D(x)$ and $D_a$ defined in \eqref{def:deltaD}--\eqref{def:Da}. Recall also that for a starting point $x\in D\setminus D_{R/2}$ (so that $\delta_D(x)\le R/2$), $z_x$ denotes the unique point in $\partial D$ such that $|x-z_x|=\delta_{D}(x)$, and $\mathbf{n}_{z_{x}}=\frac{z_{x}-x}{|z_x-x|}$ denotes the outward unit normal vector to $\partial D$ at the point $z_x$. 
With the \textit{interior $R$-ball condition} in \eqref{def:interior-exterior} in mind,
let $H_{x}$ denote the unique half-space containing the interior $R$-ball at the point $z_x$ whose normal vector is $\textbf{n}_{z_{x}}$.
Combining Lemmas \ref{lemma:half-space} and \ref{lemma:ball} gives the following statement.

\begin{lemma}\label{lemma:cancellation}
Let $D$ and $X$ be as in Theorem \ref{thm:main1}, where the $C^{1,1}$ characteristics of $D$ are given by $(R,\Lambda)$. 
Let $A\in(0,R/2]$ be the constant appearing in Lemma \ref{lemma:half-space}.
Then for any $a\in(0, A ]$,
\[
\lim_{t\to0}\psi^{-1}(1/t)\int_{D\setminus D_{a}}\P_{x}(\tau^{X}_{D}\leq t<\tau^{X}_{H_x})\ud x=0.
\]
\end{lemma}
\pf
With Lemmas \ref{lemma:half-space} and \ref{lemma:ball} at hand, the proof is similar to \cite[Lemma 3.4]{PS22}.
By \eqref{eqn:vdBD89},
$|\partial D_{a}|\leq 2^{d-1}|\partial D|$ for all $a\in(0,A]$. 
Hence,  
\begin{align*}
\int_{D\setminus D_{a}}\P_{x}(\tau^{X}_{D}\leq t<\tau^{X}_{H_x})\ud x
&\leq 2^{d-1}|\partial D|\int_{0}^{a}\P_{(\tilde{0}, u)}(\tau^{X}_{B((\tilde{0},R),R)}  \leq t<\tau^{X}_{H}   )\ud u\\
&=2^{d-1}|\partial D|\left( \int_{0}^{a}\P_{(\tilde{0}, u)}( \tau^{X}_{B((\tilde{0},R),R)} \leq t)\ud u-\int_{0}^{a}\P_{(\tilde{0}, u)}( \tau^{X}_{H} \leq t)\ud u\right).
\end{align*}
Applying Lemmas \ref{lemma:half-space} and \ref{lemma:ball} gives the desired conclusion.  
\qed

Now we are ready to derive the upper bound for Theorem \ref{thm:main1}. \bigskip

\noindent
\textbf{Derivation of the upper bound for Theorem \ref{thm:main1}.}
It follows from \eqref{eqn:vdBD89} that for any $\eps>0$, there exists $a=a(\eps)\in (0,A]$ such that 
\begin{align}\label{D-and-Du}
	|\partial D| - \eps < |\partial D_{u}|<|\partial D| +\eps \ \ \textrm{for all} \ u\leq a.
\end{align}
For this particular $a$, note that
\[
|D|-Q_{D}^{X}(t)=\int_{D_{a}}\P_{x}(\tau_{D}^{X}\leq t)\ud x +\int_{D\setminus D_{a}}\P_{x}(\tau_{D}^{X}\leq t)\ud x.
\]
Since $\psi^{-1}\in \mathcal{R}_{1/\alpha}(\infty)$ due to \cite[Theorems 1.5.3 and 1.5.12]{BGT}, it follows from Lemma \ref{lemma:inside} that 
\[
\lim_{t\to 0}\psi^{-1}(1/t)\int_{D_{a}}\P_{x}(\tau_{D}^{X}\leq t)\ud x=0.
\]
On the other hand, for any $x\in D$ with $\delta_{D}(x)\le a$ (in other words, $x\in D\setminus D_a$), 
\[
\{\tau_{D}^{X}\leq t\}\subset \{\tau^{X}_{H_x}\leq t\} \cup \{\tau_{D}^{X}\leq t <\tau_{H_x}^{X}\},
\]
where $H_x$ is the one defined right above the statement of Lemma \ref{lemma:cancellation}.
Therefore,
\begin{align}\label{inequality_split_1}
\int_{D\setminus D_{a}}\P_{x}(\tau_{D}^{X}\leq t)\ud x \leq 
\int_{D\setminus D_{a}}\P_{x}(\tau_{H_x}^{X}\leq t)\ud x
+\int_{D\setminus D_{a}}\P_{x}(\tau_{D}^{X}\leq t <\tau_{H_x}^{X})\ud x.
\end{align} 
In terms of the first integral on the right hand side, 
\[
	\int_{D\setminus D_{a}}\P_{x}(\tau_{H_x}^{X}\leq t)\ud x
	=\int_{0}^{a}|\partial D_{u}|\P_{(\tilde{0}, u)}(\tau_{H}^{X}\leq t)\ud u, 
\]
so \eqref{D-and-Du} and Lemma \ref{lemma:half-space} together gives
\begin{align}\label{inequality_split_2}
\lim_{t\to 0}\psi^{-1}(1/t)\int_{D\setminus D_{a}}\P_{x}(\tau_{H_x}^{X}\leq t)\ud x=|\partial D|\E[\overline{Y}^{(\alpha)}_{1}].
\end{align}
Combining \eqref{inequality_split_1}, \eqref{inequality_split_2}, and Lemma \ref{lemma:cancellation} yields the upper bound
\[
\limsup_{t\to0}\psi^{-1}(1/t)(|D|-Q_{D}^{X}(t))\leq |\partial D|\E[\overline{Y}^{(\alpha)}_{1}],
\]
as desired.
\qed

Now we shift our attention to the lower bound for Theorem \ref{thm:main1}, the derivation of which relies on the
\textit{exterior $R$-ball condition} in \eqref{def:interior-exterior}
for the $C^{1,1}$ open set $D$. 
In the lemma below, we consider the first exit time $\tau^{X}_{B((\tilde{0},-R),R)^{c}}$ from the \textit{complement} of the ball $B((\tilde{0},-R),R)$ of radius $R$ centered at the point $(\tilde{0},-R)$, which is located in the lower half-space 
\[
	H'=\{x=(x_1,\ldots, x_{d})\in \R^d : x_{d}< 0\}.
\]

\begin{lemma}\label{lemma:outer ball}
Let $D$ and $X$ be as in Theorem \ref{thm:main1}, where the $C^{1,1}$ characteristics of $D$ are given by $(R,\Lambda)$. 
Let $A\in(0,R/2]$ be the constant appearing in Lemma \ref{lemma:half-space}.
Then for any $a\in(0, A ]$,
\[
\lim_{t\to 0}\psi^{-1}(1/t)\int_{0}^{a}\P_{(\tilde{0},u)}(\tau^{X}_{B((\tilde{0},-R),R)^{c}}\leq t)\ud u=\E[\overline{Y}^{(\alpha)}_{1}].
\]
\end{lemma}
\pf
The proof is similar to the proof of Lemma \ref{lemma:ball}, but we provide the details for the reader's convenience. 
Fix $a\in (0,A]$.
Since $H\subset B((\tilde{0},-R),R)^{c}$, it follows that $\P_{(\tilde{0},u)}(\tau^{X}_{B((\tilde{0},-R),R)^{c}}\le t)\leq \P_{(\tilde{0},u)}(\tau_{H}^{X}\le t)$. This, together with Lemma \ref{lemma:half-space}, yields the upper bound
\[
\limsup_{t\to 0}\psi^{-1}(1/t)\int_{0}^{a}\P_{(\tilde{0},u)}(\tau^{X}_{B((\tilde{0},-R),R)^c}\leq t)\ud u\leq \E[\overline{Y}^{(\alpha)}_{1}].
\]

To derive the lower bound, let $\eps>0$ and recall that $\E[\overline{Y}_{1}^{(\alpha)}]=\int_{0}^{\infty}\P(\overline{Y}_{1}^{(\alpha)}\geq x)\ud x<\infty$ for $\alpha\in (1,2]$. There exists an integer $M=M(\eps)>0$ such that
\begin{align}\label{eqn:ex ball lb3}
\int_{0}^{M}\P(\overline{Y}_{1}^{(\alpha)}\geq x)\ud x>\E[\overline{Y}_{1}^{(\alpha)}]-\eps.
\end{align}
As in the proof of Lemma \ref{lemma:half-space}, one can verify that the law of $\tau_{B((\tilde{0},-R),R)^{c}}^{X}$ under $\P_{x}$ is equal to the law of $t\tau_{B((\tilde{0},-\psi^{-1}(1/t)R),\psi^{-1}(1/t)R)^{c}}^{\widetilde{X}^{(t)}}$ under $\P_{\psi^{-1}(1/t)x}$, where $\widetilde{X}^{(t)}$ is the scaled process defined in \eqref{def:Xtilde}.
Hence, by the change of variables $v=\psi^{-1}(1/t)u$,  
\begin{align}\label{eqn:ex ball lb1}
\int_{0}^{a}\P_{(\tilde{0},u)}(\tau^{X}_{B((\tilde{0},-R),R)^{c}}\leq t)\ud u\nn
&=\int_{0}^{a}\P_{(\tilde{0},\psi^{-1}(1/t)u)}(\tau_{B((\tilde{0},-\psi^{-1}(1/t)R),\psi^{-1}(1/t)R)^{c}}^{\widetilde{X}^{(t)}} \leq 1)\ud u\nn\\
&=\psi^{-1}(1/t)^{-1}\int_{0}^{a\psi^{-1}(1/t)/2}\P_{(\tilde{0},v)}(\tau_{B((\tilde{0},-\psi^{-1}(1/t)R),\psi^{-1}(1/t)R)^{c}}^{\widetilde{X}^{(t)}} \leq 1)\ud v\nn\\
&\geq \psi^{-1}(1/t)^{-1}\int_{0}^{M}\P_{(\tilde{0},v)}(\tau_{B((\tilde{0},-\psi^{-1}(1/t)R),\psi^{-1}(1/t)R)^c}^{\widetilde{X}^{(t)}} \leq 1)\ud v
\end{align}
for all sufficiently small $t>0$ (since $\psi^{-1}(1/t)\to\infty$ as $t\to 0$). 

Now we construct an increasing sequence of domains $E(n)$ by
\[
E(n)=\left\{x=(\tilde{x}, x_{d})\in \R^{d-1}\times \R: |\tilde{x}|<n,\, -n<x_{d}< 0,\, \cos^{-1}(-\vec{\bf{n}}\cdot \frac{x}{|x|} )<\frac{\pi}{2}-\frac{1}{n} \right\}, \quad -\vec{\bf{n}}=(\tilde{0},-1).
\]
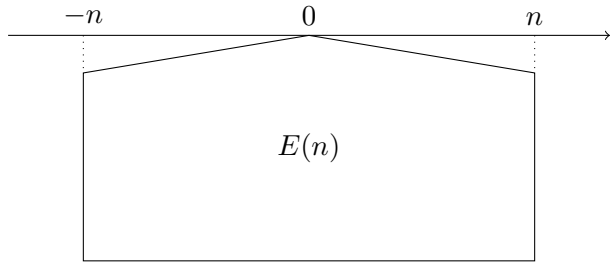
\begin{center}
\begin{tikzpicture}
\draw [->] (-4,0)--(4,0);
\draw (0,0)--(-3,-0.5)--(-3,-3)--(3,-3)--(3,-0.5)--(0,0);
\draw [dotted] (3,0)--(3,-0.5);
\draw [dotted] (-3,0)--(-3,-0.5);
\node [above] at (3,0) {$n$};
\node [above] at (-3,0) {$-n$};
\node [above] at (0,0) {0};
\node at (0,-1.5) {$E(n)$};
\end{tikzpicture}
\end{center}
Since $E(n)$ increases to the lower half-space $H'=\{(x_1,\ldots, x_{d}) : x_{d}< 0\}$ as $n\to\infty$, it follows from the bounded convergence theorem that 
\[
\lim_{n\to\infty}\int_{0}^{M}\P_{(\tilde{0},v)}(\tau_{E(n)^c}^{X^{(\alpha)}}\leq 1)\ud v=\int_{0}^{M}\P_{(\tilde{0},v)}(\tau_{(H')^{c}}^{X^{(\alpha)}}\leq 1)\ud v
=\int_{0}^{M}\P_{(\tilde{0},v)}(\tau_{\overline{H}}^{X^{(\alpha)}}\leq 1)\ud v.
\]
Since the latter sequence converges increasingly, we may take an integer $N=N(\eps)$ such that 
\beq\label{eqn:ex ball lb2}
\int_{0}^{M}\P_{(\tilde{0},v)}(\tau_{E(N)^{c}}^{X^{(\alpha)}}\leq 1)\ud v>\int_{0}^{M}\P_{(\tilde{0},v)}(\tau_{\overline{H}}^{X^{(\alpha)}}\leq 1)\ud v - \eps.
\eeq
Since $B((\tilde{0},-\psi^{-1}(1/t)R),\psi^{-1}(1/t)R)$ increases to $H'$ as $t\to 0$, we can take $t_{0}=t_{0}(N)>0$ such that for any $0<t\leq t_{0}$,
\[
E(N)  \subset B((\tilde{0},-\psi^{-1}(1/t)R),\psi^{-1}(1/t)R), 
\]
which implies that 
\begin{align}\label{eqn:ex ball lb4}
\P_{(\tilde{0},v)}(\tau_{B((\tilde{0},-\psi^{-1}(1/t)R),\psi^{-1}(1/t)R)^{c}}^{\widetilde{X}^{(t)}} \leq 1)
\geq \P_{(\tilde{0},v)}(\tau_{E(N)^{c}}^{\widetilde{X}^{(t)}}\leq 1).
\end{align}

Combining \eqref{eqn:ex ball lb1}--\eqref{eqn:ex ball lb4} and using Fatou's lemma and the weak convergence of $\widetilde{X}^{(t)}$ to $X^{(\alpha)}$, we obtain
\begin{align*}
\liminf_{t\to0}\psi^{-1}(1/t)\int_{0}^{a}\P_{(\tilde{0},u)}(\tau^{X}_{B((\tilde{0},-R),R)^{c}}\leq t)\ud u
&\geq \liminf_{t\to0}\int_{0}^{M}\P_{(\tilde{0},v)}(\tau_{E(N)^{c}}^{\widetilde{X}^{(t)}} \leq 1)\ud v\\
&\ge \int_{0}^{M}\P_{(\tilde{0},v)}(\tau_{E(N)^{c}}^{X^{(\alpha)}}\leq 1)\ud v\\
&> \int_{0}^{M}\P_{(\tilde{0},v)}(\tau_{\overline{H}}^{X^{(\alpha)}}\leq 1)\ud v - \eps\\
&= \int_{0}^{M}\P(\overline{Y}^{(\alpha)}_{1}> v)\ud v-\eps
> \E[\overline{Y}_{1}^{(\alpha)}]-2\eps,
\end{align*}
where the last line follows from \eqref{eqn:ex ball lb3}.
Since $\eps>0$ is arbitrary, we conclude that 
\[
\liminf_{t\to0}\psi^{-1}(1/t)\int_{0}^{a}\P_{(\tilde{0},u)}(\tau^{X}_{B((\tilde{0},-R),R)^{c}}\leq t)\ud u\geq \E[\overline{Y}^{(\alpha)}_{1}],
\]
as desired. 
\qed

Recall that for a starting point $x\in D\setminus D_{R/2}$, $H_{x}$ denotes the half-space containing the interior $R$-ball at the point $z_x$ and tangent to $\partial D$. Let $\tilde{B_x}$ denote the unique exterior $R$-ball in $D^{c}$ touching the point $z_x$. (Note that $H_x\subset \tilde{B_x}^{c}$.) 
The following lemma is similar to Lemma \ref{lemma:cancellation}.

\begin{lemma}\label{lemma:cancellation2}
Let $D$ and $X$ be as in Theorem \ref{thm:main1}, where the $C^{1,1}$ characteristics of $D$ are given by $(R,\Lambda)$. 
Let $A\in(0,R/2]$ be the constant appearing in Lemma \ref{lemma:half-space}.
Then for any $a\in(0, A ]$,
\[
\lim_{t\to 0}\psi^{-1}(1/t)\int_{D\setminus D_{a}}\P_{x}(\tau_{H_x}^{X}\leq t <\tau_{\tilde{B_x}^{c}}^{X} )\ud x=0. 
\]
\end{lemma}
\pf
By \eqref{eqn:vdBD89},
$|\partial D_{a}|\leq 2^{d-1}|\partial D|$ for all $a\in(0,A]$. Hence,
\begin{align*}
\int_{D\setminus D_{a}}\P_{x}(\tau_{H_x}^{X}\leq t <\tau_{\tilde{B_x}^{c}}^{X} )\ud x
&=\int_0^{a} |\partial D_u| \P_{(\tilde{0},u)} (\tau_{H}^{X}\leq t <\tau_{B((\tilde{0},-R),R)^{c}}^{X} )\ud u\\
&\le 2^{d-1} |\partial D|\int_0^{a} \P_{(\tilde{0},u)} (\tau_{H}^{X}\leq t <\tau_{B((\tilde{0},-R),R)^{c}}^{X} )\ud u\\
&= 2^{d-1} |\partial D|\left(\int_0^{a} \P_{(\tilde{0},u)} (\tau_{H}^{X}\leq t)\ud u -\int_0^{a} \P_{(\tilde{0},u)} (\tau_{B((\tilde{0},-R),R)^{c}}^{X} \le t)\ud u\right).
\end{align*}
Applying Lemmas \ref{lemma:half-space} and \ref{lemma:outer ball} gives the desired conclusion.  
\qed

Now we are ready to derive the lower bound for Theorem \ref{thm:main1}, which extends \cite[Lemma 3.1]{PS22} for $\alpha$-stable processes with $\alpha\in (1,2]$ to any isotropic L\'evy processes with characteristic exponents $\psi\in \mathcal{R}_{\alpha}(\infty)$.
We remark that there is a minor mistake in the proof of \cite[Lemma 3.1]{PS22}. 
The authors there claimed that $\{\tau_{H_{x}}^{X^{(\alpha)}}\leq t\}\subset \{\tau_{D}^{X^{(\alpha)}}\leq t\}$, which is not true when $D$ is non-convex since the containment $D\subset H_x$ does not necessarily hold; 
we take this opportunity to fix their argument by observing the containments in \eqref{eqn:containments} below.
\bigskip

\noindent
\textbf{Derivation of the lower bound for Theorem \ref{thm:main1}.}
As we did in the derivation of the upper bound for Theorem \ref{thm:main1}, for a fixed $\eps>0$, take $a=a(\eps)\in(0, A ]$ such that \eqref{D-and-Du} holds, which results in \eqref{inequality_split_2}. 

Under $\P_x$ with $x\in D\setminus D_a$, 
\begin{align}\label{eqn:containments}
\{\tau_{H_x}^{X}\leq t\} 
&\subset \{\tau_{H_x}^{X}\leq t, \tau_{D}^{X} \leq t\} \cup \{\tau_{H_x}^{X}\leq t<\tau_{D}^{X} \}\notag\\
&\subset  \{\tau_{D}^{X} \leq t\} \cup \{\tau_{H_x}^{X}\leq t <\tau_{D}^{X} \}\notag\\
&\subset  \{\tau_{D}^{X} \leq t\} \cup \{\tau_{H_x}^{X}\leq t <\tau_{\tilde{B_x}^{c}}^{X} \},
\end{align}
where $H_x$ and $\tilde{B_x}$ are the ones defined right above the statement of Lemma \ref{lemma:cancellation2}. 
The latter implies that
\[
\P_{x}(\tau_{H_x}^{X}\leq t)\leq \P_{x}(\tau_{D}^{X} \leq t) + \P_{x}(\tau_{H_x}^{X}\leq t <\tau_{\tilde{B_x}^{c}}^{X} ).
\]
Hence, 
\begin{align*}
|D|-Q^{X}_{D}(t)
&=\int_{D}\P_{x}(\tau^{X}_{D}\leq t)\ud x
\geq \int_{D\setminus D_{a}}\P_{x}(\tau^{X}_{D}\leq t)\ud x\\
&\geq \int_{D\setminus D_{a}}\P_{x}(\tau_{H_x}^{X}\leq t)\ud x-\int_{D\setminus D_{a}}\P_{x}(\tau_{H_x}^{X}\leq t <\tau_{\tilde{B_x}^{c}}^{X} )\ud x.
\end{align*}
Combining \eqref{inequality_split_2} and Lemma \ref{lemma:cancellation2} gives the lower bound
\[
\liminf_{t\to0}\psi^{-1}(1/t)(|D|-Q_{D}^{X}(t))\geq |\partial D|\E[\overline{Y}^{(\alpha)}_{1}],
\]
as desired.
\qed

\section{Examples}\label{section:examples}
This section illustrates the fact that our results, Theorem \ref{thm:main1} and Corollary \ref{cor:main}, are applicable to a broad range of L\'evy processes.  

\begin{example}[Pure Jump Process]\label{example:pure jump}
\begin{em}
Theorem \ref{thm:main1} applies to the following pure jump processes:
\begin{enumerate}
\item Stable Process $\psi(\xi)=|\xi|^\alpha$, where $\alpha\in (1, 2)$.
\item Mixed Stable Process $\psi(\xi)=|\xi|^\alpha+|\xi|^\beta$,  where $0<\beta<\alpha$ and $\alpha\in (1,2)$.
\item Relativistic Stable Process $\psi(\xi)=(|\xi|^2+m^{2/\alpha})^{\alpha/2}-m$,  where $m>0$ and $\alpha\in (1,2)$.
\item $\psi(\xi)=|\xi|^\alpha(\log(|\xi|^2+1))^{\beta/2}$,  where $\alpha\in (1, 2)$ and $0<\beta<2-\alpha$.
\item $\psi(\xi)=|\xi|^\alpha(\log(|\xi|^2+1))^{-\beta/2}$,  where $\alpha\in (1,2)$ and $0<\beta<\alpha$.
\end{enumerate}

More generally, let $X:=\{W_{S_t}\}_{t\ge 0}$ be a subordinate Brownian motion, where $S_{t}$ is a subordinator (an increasing L\'evy process on $\R$) with Laplace exponent $\phi\in \RR_{\alpha/2}(\infty)$ with $\alpha\in (1,2]$ that is independent of a Brownian motion $W$. Then the characteristic exponent $\psi$ of $X$ takes the form $\psi(\xi)=\phi(|\xi|^2)$, and in particular, $\psi\in\RR_{\alpha}(\infty)$. Therefore, Theorem \ref{thm:main1} is applicable to a large class of subordinate Brownian motions.
For more details on subordinators, readers are referred to \cite{SSV}. 
\end{em}
\end{example}

\begin{example}[Jump Diffusion]\label{example:jump diffusion}
\begin{em}
Let $\psi_1(\xi)=|\xi|^2 + \psi(\xi)$, where $\psi$ is the characteristic exponent in Example \ref{example:pure jump}. 
This corresponds to the independent sum of a Brownian motion and a pure jump process. Clearly, $\psi_1\in \mathcal{R}_{2}(\infty)$, so Theorem \ref{thm:main1} is applicable. 
\end{em}
\end{example}

\begin{example}[Truncated Jump Process and Truncated Jump Diffusion]
\begin{em} 
Let $\nu$ be the L\'evy measure corresponding to a characteristic exponent $\psi$ appearing in Example \ref{example:pure jump}. Define a new L\'evy measure $\mu$ by
\[
\mu(\ud y):=\nu(\ud y)\cdot 1_{\{|y|\leq 1\}}.
\]
The L\'evy process with L\'evy measure $\mu$ corresponds to a truncated jump process (a process whose jump sizes are at most 1).
Since any L\'evy measure is integrable away from the origin, $\nu-\mu$ is a finite (non-negative) measure, and hence, Corollary \ref{cor:main} is applicable. 
A similar discussion applies to truncated jump diffusions constructed from jump diffusions in Example \ref{example:jump diffusion}. 
\end{em}
\end{example}

\bigskip
\noindent
{\bf Acknowledgments:}
The authors are grateful to Professor Renming Song (University of Illinois at Urbana-Champaign) for his valuable input regarding the material presented in this paper. 
They also appreciate the comments and suggestions from the anonymous referee that helped improve the exposition of the paper.
\bigskip
\noindent

\begin{singlespace}

\end{singlespace}
\end{doublespace}

\vskip 0.3truein

{\bf Kei Kobayashi}

Department of Mathematics, Fordham University, New York, NY 10023, 
USA

E-mail: \texttt{kkobayashi5@fordham.edu}

\bigskip

{\bf Hyunchul Park}

Department of Mathematics, State University of New York at New Paltz, NY 12561,
USA

E-mail: \texttt{parkh@newpaltz.edu}
\end{document}